\documentclass[5p]{elsarticle}

\usepackage{amssymb}
\usepackage{amsmath}
\usepackage{color}
\usepackage[normalem]{ulem}
\usepackage{tikz}
\usetikzlibrary{patterns}
\journal{Journal of Computational Physics}

\begin{document}

\begin{frontmatter}




\title{A GPU cluster optimized multigrid scheme for computing unsteady incompressible fluid flow}

\author{Gy\"orgy Tegze\corref{cor1}}
\ead{tegze.gyorgy@wigner.mta.hu}
\cortext[cor1]{Corresponding author, Tel: (+36) 1 392 2222 ext. 3839; Fax: (+36) 1 392 2219}

\author{Gyula I. T\'oth\corref{cor2}}

\address{Institute for Solid State Physics and Optics, Wigner Research Centre for Physics
P.O. Box 49, H-1525 Budapest, Hungary}

\begin{abstract}
A multigrid scheme is proposed \textcolor{black}{for the pressure equation of the incompressible unsteady fluid flow equations, allowing} efficient implementation on \textcolor{black}{clusters of} modern CPUs, many integrated core devices (MICs), and graphics processing units (GPUs).
\textcolor{black}{It is shown that the total number of the synchronization events can be significantly reduced} when a deep, 2h grid hierarchy is replaced with a two-level scheme using 16h-32h restriction, fitting to the the width of the SIMD engine of modern CPUs and GPUs.
In addition, optimal memory transfer is also ensured, since no strided memory access is required.
\textcolor{black}{We report increasing arithmetic intensity of the smoothing steps when compared to the conventional additive correction multigrid (ACM), however it is counterbalanced in runtime by the decreasing number } of the expensive restriction steps.
\textcolor{black}{A systematic construction methodology for the coarse grid stencil is also presented that helps in moderating the excess arithmetic intensity associated with the aggressive coarsening.}
Our higher order interpolated stencil improves convergence rate via minimizing spurious interference between the coarse and the fine scale solutions.
The method is demonstrated on solving the pressure equation for 2D incompressible fluid flow:
The benchmark setups cover shear driven laminar flow in cavity, and direct numerical simulation (DNS) of a turbulent jet.
\textcolor{black}{We have compared our scheme to the ACM in terms of the arithmetic intensity of the iterations and the number of the synchronization calls required.
Also the strong scaling is plotted for our scheme} when using a hybrid OpenCl/MPI based parallelization.
%structured grid
\end{abstract}


\end{frontmatter}

% \linenumbers

\section{Introduction}

Numerical studies of incompressible fluid flow are important in both academic research and engineering applications.
The incompressibility constraint and constant density represent a good approximation when fluid flow velocity is significantly smaller than the speed of sound in the media.
For subsonic unsteady flows this approach usually allows larger time-steps, because, contrary to the compressible fluid models, density waves travelling at the speed of sound do not need to be resolved.
However, solving for the pressure instead of density changes the  continuity equation from parabolic to elliptic type.
The major computational challenge in simulating incompressible unsteady fluid flow is to \textcolor{black}{develop an effective, parallel pressure equation solver} \cite{Kwak_2005_CANDF}.

The demand for more detailed descriptions of convection phenomena, such as microscale description of multiphase systems, large eddy simulation (LES) \cite{Smagorinsky_1963_MWR,Deardorff_1970_JFM} or direct numerical simulation (DNS) \cite{Orszag_1970_JFM}, drives evolution of mathematical and numerical tools along with the evolving computer architectures.
Recently, the share of accelerator cards (i.e. GPUs), and multi-core coprocessors is rapidly increasing in the supercomputing scene.
Alongside, programming paradigms are changing to fully exploit the fine grained parallelism available on these hardwares.
Writing GPU optimized codes for solving parabolic PDEs using explicit time marching is straightforward \cite{mick_2009_GPGPU}, while computationally efficient treatment of elliptic problems such as the Poisson equation \textcolor{black}{together with} incompressibility requires more complex methods.
The difficulty lies in the multi-level parallelism of the current supercomputing hardwares.
While efficient methods exist for large scale elliptic problems on distributed memory hardware, it is not straightforward to optimize these algorithms for \textcolor{black}{both fine-grain and coarse-grain parallelism coincidentally \cite{Muller_2014_arxiv,Tarjan_1996_IJPP}.  
Although multiple compute speed is available on recent accelerator cards such as GPUs, they are accompanied with significant latency when copying data between such devices.
A simple synchronization between cards involves a device to host, a host to host, and a host to device memory transfer, which at least \textit{triples synchronization latency.}
It is not uncommon, that the ratio of the computational time and the memory transfer latency increases by even 2 orders of magnitude compared to CPU clusters.
Therefore, in the present work we pay special attention to optimize our scheme for the case when latency limits the parallel efficiency of the computations.
}

Multigrid (MG) methods \textcolor{black}{ \cite{Ghia_1982_JCP} are developed for to} accelerate the iterative solution of large algebraic equations, such as the discretized pressure equation.
The number of iterations required to reach the convergence criteria \textcolor{black}{is reduced by using a} multi-resolution discretization that \textcolor{black}{accelerates the relaxation of long-wavelength components of the residual.
Compared to simple iterative methods, MG has a low memory footprint \cite{Muller_2014_arxiv}, which is a significant benefit on accelerator devices having limited onboard memory.}
The MG methods traditionally use a deep hierarchy of discretization, where each grid is 2h coarser than the previous level.
The scheme development for this multi-resolution hierarchy needs special attention, since parallel execution of the iterations on fine meshes are usually bandwidth limited, while the solution of the coarser level discrete equations can be latency limited.
Besides, it is not straightforward to achieve optimal fine grain parallelism when using these schemes.
It is difficult to optimize memory load performance, cache usage and SIMD utilization at the same time, even if the programming model allows control over these procedures.
The ease of programming is also an issue when evaluating a numerical procedure.
Here, not only field variables, but also boundary conditions must be consecutively transfered through the grid hierarchy ensuring proper padding to keep data alignment.
%Compared to other programming languages, load balancing in streaming languages requires different approach, since it is not possible to launch and terminate "threads" inside compute kernels.
%Although the conventional 2h hierarchy of the multigrid scheme does not fit to streaming programming languages, it can be implemented on GPU architectures
Despite challenges, various multigrid schemes have been successfully implemented on GPUs \cite{Bolz_ACM_TOG,Bell_2012_SIAM}, and GPU clusters \cite{Cohen_2013_PE,Goddeke_2008_IJCSE}.

Here we present an aggressive coarsening strategy, that results in a simple but effective two-level multigrid scheme.
Our scheme bypasses some of the programming difficulties arising when implementing to many-core hardware architectures, and decrease the memory footprint of the solver.
Our method is compared to the Additive Correction Multigrid (ACM) method for unsteady incompressible fluid flows.
Simulation of a laminar flow in a shear driven cavity, and a direct numerical simulation of a turbulent jet have been used to benchmark against the ACM scheme.
Finally we present strong scaling figures for a channel-flow setup, that mimics 2D turbulence in soap film experiments \cite{Kellay_1995_PRL}.

\section{The fluid model and its discretization}
\subsection{The equations of the incompressible fluid}
Three assumptions on fluid behavior yields the Navier-Stokes equation: (i) The viscous dissipation is a linear function of the strain rates, (ii) the fluid is isotropic (rotational invariance), and (iii) for a fluid at rest hydrostatic pressure applies.
Assuming incompressibility and constant density, the Navier-Stokes equation and the continuity of mass read as follows:
\begin{eqnarray}
	\frac{\partial \mathbf{v}} {\partial t} &=& \mathbf{v} \cdot ({\nabla} \otimes \mathbf{v})+ \eta\Delta\otimes\mathbf{v} +\nabla p \label{eq_NS}\\
	0 &=& -{\nabla} \cdot \mathbf{v}\label{eq_cont}
\end{eqnarray}
where equation (\ref{eq_NS}) is solved for the fluid velocity $\mathbf v$, and equation (\ref{eq_cont}) is solved for the pressure $p$ when substituting $\mathbf v$ [see Eq. (4)].
$\eta$ is the kinematic viscosity.

\subsection{The time stepping scheme}
%Chorin's projection method \cite{Chorin_1968_MATHCOMP} was chosen to solve the governing equations.
\textcolor{black}{A simple variant of Chorin's projection method \cite{Chorin_1968_MATHCOMP}, the Goda's incremental pressure correction method \cite{Goda_1979_JCP} has been chosen to solve the governing equations.}
We predict velocity $\mathbf v^*$ for the next time-step using a simple first order time integration of the explicitly known terms.
The pressure can be decomposed as: $p^{t+1}=p^t+\delta p$.
Using $p^t$ the predicted velocity reads as:
\begin{equation}
	\mathbf v^{*} = \mathbf v^{t} + \Delta t [\mathbf{v}^t \cdot ({\nabla} \otimes \mathbf{v}^t)+ \eta\Delta\otimes\mathbf{v}^t] +\nabla p^t .
\end{equation}
Substituting $\mathbf v^{t+1} = \mathbf v^* + \nabla \delta p$ into Eq. (\ref{eq_cont}) yields the following Poison equation:
\begin{equation}
	0 = \nabla \cdot \mathbf{v}^{*} + \nabla^2 \delta p ,
\end{equation}
which has to be solved for the pressure change $\delta p$.
Finally $\mathbf v^*$ is corrected with the pressure-change term $\nabla \delta p$ to obtain the divergence free $\mathbf{v}^{t+1}$.

\subsection{The spatial discretization}
A staggered grid arrangement \cite{Harlow_1965_POF} has been chosen for the velocity to avoid odd-even decoupling and the resulting spurious checkerboard patterns in the solution.
Finite differences of second order accuracy were used to discretize the Navier-Stokes equation \cite{E_Weinan_2001_MATHCOMP}, and a five point Laplacian stencil was used to discretize the pressure equation.
Alternatively, the discretization can be derived from finite volumes (FV) that gives a clear explanation of constructing flux preserving conservative stencils.

\section{Solution of the pressure equation}
\subsection{Simple methods}
While the explicit time marching scheme for the Navier-Stokes equation can be solved line by line, discretizing the continuity equation leads to a large system of linear equations.
Direct solution of such systems is computationally too expensive for most practical problems, and also difficult to perform in a parallel manner.
The discretization error is usually larger than the accuracy of the computer arithmetic, thus an approximate solution of the discrete problem is satisfactory, and an iterative solution can be applied.
In parallel environment, the Gauss-Sedel (GS) iteration is commonly applied.
Concurrent computing of the well known 5 point stencil is feasible when using red-black reordering \cite{Kuo_1990_SIAM_SSC}, while  multicolor ordering \cite{Adams_1988_SIAM_JNA} can be used for higher order stencils.
%Using red-black reordering \cite{Kuo_1990_SIAM_SSC} the well known 5 point stencil can be executed in parallely, while  multicolor ordering \cite{Adams_1988_SIAM_JNA} can be used for higher order stencils.

\subsection{The multigrid method}
The GS iterations are very effective in removing the high frequency components of the error, but the low frequency components decay with a very low rate.
Convergence requires $O(N^2)$ iterations, where $N$ is the linear size.
Multigrid methods are commonly used to accelerate convergence: coarser discretizations are used to eliminate the lower frequency components of the residual.
The MG technique consists of a sequence of smoothing, restriction, and prolongation operators.
Smoothing steps are the actual iterations using a solver (i.e. Gauss-Seidel) on a grid level.
Restriction is a downsampling the residual to a coarser discretization, and prolongation is an interpolation from a coarser to finer discretization.
A typical multigrid cycle starts with smoothing of the residual error on the finest level, then the approximate solution is transfered to a coarser level (restriction), and smoothened again.
The procedure is repeated until the coarsest level is reached, then the way back.
The optimal depth of the grid hierarchy depends on the actual problem, the  discretized equations must be solved in a reasonable time on the coarsest level.
A direct solver is often used at the coarsest grid.

When using a geometric multigrid (GMG) scheme, the original PDE is re-discretized on each grid, data structures must be constructed on each level.
Note that the solutions on coarse grids does not necessarily approximate the solution on the fine grid, but it approximates the original problem.
Possible inconsistency between fine and coarse grid solutions may results in extra iteration cycles before reaching convergence.

When using an algebraic multigrid (AMG) scheme, agglomeration of the coefficient matrix is performed.
In the case of the widely used AMG procedure, the Additive Correction Multigrid  scheme \cite{Hutchinson_1986_NHT}, the coarse grid equations are obtained by the summation of the finer grid equations.
The results obtained on the coarser grid are simply added to the finer grid solution, no interpolation or extrapolation operators are required in this procedure.
The real advantage of this method is that on each level it approximates the solution of the discrete problem on the finest grid.
Also note that conserving property of the spatial discretization scheme is satisfied at all discretization level.

\subsection{\textcolor{black}{The construction of an interpolated stencil for the coarse mesh.}}
We extend the basic idea of the ACM scheme in our MG procedure.
It is often useful to interpret the discretization of the conservation equations (either momentum or mass) in the finite volume manner.
When constructing a FV scheme, both sides of the PDE are integrated over a control volume (CV).
Applying the Gauss theorem on the homogeneous part of the Poisson equation results in a surface integral of the fluxes over the CV faces.
Thus simplest discretization is obtained, when the inhomogeneous part is assumed to be the CV average, and the surface normals of the pressure gradients are assumed to be constant at all CV faces.
This recovers the well-known five point Laplacian stencil.
The overall accuracy of the respective FV scheme depends both on the accuracy of calculating the fluxes on the CV face (using finite differences), and the accuracy of the surface integral \cite{Peric_CFD}.

\begin{figure}[ht!]
\begin{center}
\definecolor{cccccc}{rgb}{0.8,0.8,0.8}
\begin{tikzpicture}[line cap=round,x=1.0cm,y=1.0cm]

%compute cells
\foreach \j in {0, 0.5, ..., 6.5}
\foreach \i in {0, 0.5, ..., 7.5} {
	\draw[color=orange] (\i,\j) rectangle (\i+0.5,\j+0.5);
}
%north boundary padding
\foreach \j in {0, 0.5, ..., 1}
\foreach \i in {0, 0.5, ..., 8} {
	\filldraw[draw=white,fill=orange] (\i,7+\j) rectangle (\i+0.5,7.5+\j);
}
%east boundary padding
\foreach \i in {0, 0.5, ..., 8} {
	\filldraw[draw=white,fill=orange] (8,\i) rectangle (8.5,\i+0.5);
}
%vertical tile boundaries
\foreach \i in {0, 4, ..., 8} {
	\draw [line width=1.5pt] (\i+0.5,0)-- (\i+0.5,8.5);
}
%horisontal tile boundaries
\foreach \i in {0, 4, ..., 8} {
	\draw [line width=1.5pt] (0,\i+0.5)-- (8.5,\i+0.5);
}

%coarse grained cell
\draw [line width=0.5pt,pattern=north west lines, pattern color=orange] (0.5,0.5) rectangle (4.5,4.5); 
%interpolation domain
\draw [dashed,line width=1.5pt] (2.5,2.5) rectangle (2.5+3.75,2.5+3.25); 

% SW interp center
\definecolor{mycolor}{rgb}{0,0,0};
\draw [fill=mycolor] (2.5,2.5) circle [radius=.2];
\draw [font=\small] (2.5,2.5) node[anchor=north] {$Q_{11}=P_C$};

% SW interp center
\definecolor{mycolor}{rgb}{0.1,0.1,0.1};
\draw [fill=mycolor] (2.5+3.75,2.5) circle [radius=.2];
\draw [font=\small] (2.5+3.75,2.5) node[anchor=north] {$Q_{21}=P_E$};

\definecolor{mycolor}{rgb}{0.5,0.5,0.5};
\draw [fill=mycolor] (2.5,2.5+3.25) circle [radius=.2];
\draw [font=\small] (2.5,2.5+3.25) node[anchor=north] {$Q_{12}=P_N$};

\definecolor{mycolor}{rgb}{0.99,0.99,0.99};
\draw [fill=mycolor] (2.5+3.75,2.5+3.25) circle [radius=.2];
\draw [font=\small] (2.5+3.75,2.5+3.25) node[anchor=north] {$Q_{22}=P_{NE}$};

%interpolated values
\foreach \j in {1, 2, ..., 4} {
\foreach \i in {4, 5} {
	%delta x = 7.5
	%delta y = 6.5
	%Q11 = 0.0; Q21 = 0.9; Q12 = 0.8; Q22 = 0.95
	\pgfmathparse{1.0/(7.5*6.5) * ( 0.1*(\i-0.5)*(6.5-\j+0.5) + 0.5*(7.5-\i+0.5)*(\j-0.5) + 0.99*(\i-0.5)*(\j-0.5) )}
	\definecolor{currentcolor}{rgb}{\pgfmathresult, \pgfmathresult, \pgfmathresult};
	\filldraw[draw=white,fill=currentcolor] (2+0.5*\i,2+0.5*\j) rectangle (2.5+0.5*\i,2.5+0.5*\j);
}
}

%interpolated values
\foreach \j in {4,5} {
\foreach \i in {1, 2, ..., 4} {
	%delta x = 7.5
	%delta y = 6.5
	%Q11 = 0.0; Q21 = 0.9; Q12 = 0.8; Q22 = 0.95
	\pgfmathparse{1.0/(7.5*6.5) * ( 0.1*(\i-0.5)*(6.5-\j+0.5) + 0.5*(7.5-\i+0.5)*(\j-0.5) + 0.99*(\i-0.5)*(\j-0.5) )}
	\definecolor{currentcolor}{rgb}{\pgfmathresult, \pgfmathresult, \pgfmathresult};
	\filldraw[draw=white,fill=currentcolor] (2+0.5*\i,2+0.5*\j) rectangle (2.5+0.5*\i,2.5+0.5*\j);
}
}

\foreach \i in {1, 2, ..., 4} {
	\draw [->,,>=stealth,white,thick] ((2.25+0.5*\i,4.25) -- +(0,0.5);
}
\foreach \j in {1, 2, ..., 4} {
	\draw [->,,>=stealth,white,thick] (4.25,2+0.5*\j+0.25) -- +(0.5,0);
}

\end{tikzpicture}

\end{center}

\caption{\label{fig:mg_grid}
Multigrid layout: White rectangles show fine grid-cells, and orange coloring indicate padding to help coalesced read and write.
Coarse grid is drawn by solid black lines, some cell centers ($P_C$,$P_E$,$P_N$ and $P_{NE}$) are denoted with shaded circles.
The actual center node for the coarse grid is hatched.
The interpolation area for the north-east corner of the center node is indicated using dashed lines, and fine grid cells whose interpolated values used to construct the coarse grid discretized equations are shaded.
White arrows show pressure gradients that are summed for the coarse grid discretization.}

\end{figure}
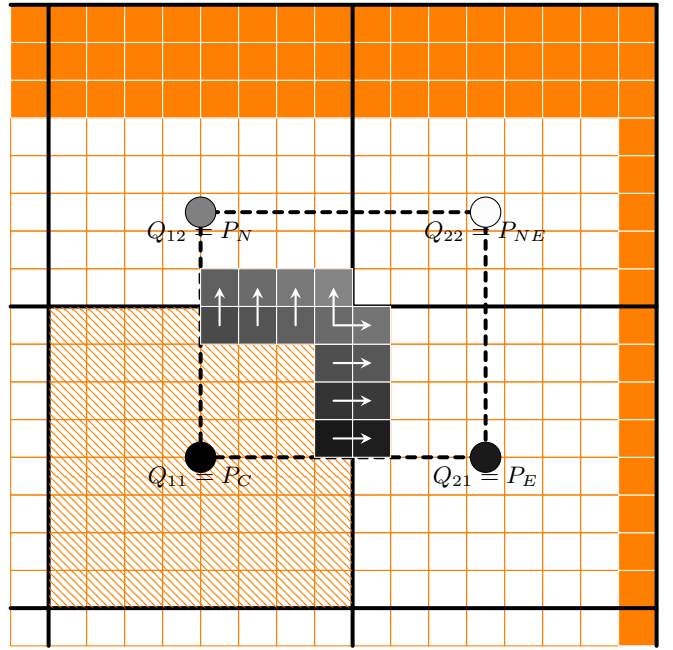

After restricting the fine grid values to a coarser grid cell (hatched area in Fig. \ref{fig:mg_grid}.), the same finite volume discretization procedure applies when the traditional geometric multigrid (GMG) method is used.
However, one can use many values that are already pre-calculated on the fine mesh to make the discretization more accurate and easier to compute.
For example, instead of downsampling $\mathbf v^*$ to the coarse mesh and calculate its divergence using the corresponding finite volume stencil, one can sum $\nabla \mathbf v^*$ already calculated for the fine grid discretization.
In the latter case the integration for the coarse CV will be the same accuracy as for the fine grid.
The procedure is identical to the agglomeration of the coefficient matrix in the ACM method.
Naturally, the integration (summation) can be applied for larger volumes too (a $8h$ restriction is shown in Fig. 1).
%2h restriction is regularly applied in the ACM method, however it can be done on larger volumes (8h restriction is shown in Fig. \ref{fig:mg_grid}).
Moreover one can perform similar sum for the fluxes derived from applying the Gauss theorem for homogeneous part of the PDE.
It is clear that fluxes not crossing the coarse grid cell faces cancel each other in the volumetric integral (summation), therefore, we have to account for the discretization only at the cell faces.
%It is clear that for a conservative scheme, the fluxes that do not cross the coarse grid cells are cancelling each other.
%We just have to account for the discretization at the cell face.
When applying to domains of arbitrary size, an interpolation rule must be constructed  to define fine grid pressure values $p_{i,j}$ near the boundaries.
Thus the finite volume scheme for the coarse grid is constructed by summing fluxes calculated from the $p_{i,j}$ interpolated values.
At this point the method differs from the ACM scheme, where a simple summation is performed that corresponds to assuming constant $p_{i,j}$ in the entire coarse CV.
%To maintain the conservative character, the neighboring coarse cells should use the same interpolated values.
The neighboring coarse cells must use the same interpolated values to maintain the conservative character.
This can be easily ensured, if these cells use a common domain to interpolate over.
A natural choice is to interpolate between four neighboring cell centers (dashed rectangle in Fig. \ref{fig:mg_grid}.).
We have used a simple bilinear interpolation to calculate $p_{i,j}$ as the function of pressure values at the coarse cell centers $P_{i,j}$.
Naturally, the same interpolation was used to prolongate the pressure values to the fine grid.
If the origin is ($Q_{11}$) the bilinear interpolation reads as:

%\begin{subequations}
%\begin{align}
\begin{eqnarray}
\label{eq_bilinear interp}
p(x,y) &=& \frac{1}{\Delta x^{W|E} \Delta y^{S|N}} [ \\
					&&Q_{11} (\Delta x^{W|E}-x)(\Delta y^{S|N} -y) + \nonumber\\
                                        &&Q_{21} x(\Delta y^{S|N} -y) + \nonumber\\
                                        &&Q_{12} (\Delta x^{W|E}-x) y +\nonumber \\
                                        &&Q_{22} x y \nonumber] ,
\end{eqnarray}
%\end{align}
%\end{subequations}
where $\Delta x^{W|E}$ and $\Delta y^{S|N}$ are the horizontal, and the vertical size of the interpolation rectangle, while $x$ and $y$ are the coordinates of the fine grid control volume centers.
The superscripts indicate that the size of the interpolation rectangle can differ on the west-east and the south north faces, because of padding, whereas the coarse mesh can be nonuniform. 
$Q_{I,J}$ are pressure values at the corners.
Now, the fluxes across fine CV faces $S_i$ and $S_j$ (white arrows in Fig. \ref{fig:mg_grid}.) can be computed as follows:
%\begin{subequations}
%\begin{align}
\begin{eqnarray}
\label{eq_surf_xderiv_2}
\frac{\partial p(x,y)}{\partial x}\Bigg\vert_{S_j} &=& \frac{1}{\Delta x^{W|E} \Delta y^{S|N}} [\\
        &&Q_{11} \{-1\} (\Delta y^{S|N} -y_j) +  \nonumber\\
        &&Q_{21} \{+1\}(\Delta y^{S|N} -y_j) +  \nonumber\\
        &&Q_{12} \{-1\} y_j + \nonumber \\
        &&Q_{22} \{+1\} y_j ] \nonumber
\end{eqnarray}
%\end{align}
%\end{subequations}
and
%\begin{subequations}
%\begin{align}
\begin{eqnarray}
\label{eq_surf_yderiv_1}
\frac{\partial p(x,y)}{\partial y}\Bigg\vert_{S_i} &=& \frac{1}{\Delta x^{W|E} \Delta y^{S|N}} [\\
        &&Q_{11} \{-1\} (\Delta x^{W|E} -x_i) + \nonumber \\
        &&Q_{21} \{-1\} x_i + \nonumber\\
        &&Q_{12} \{+1\}(\Delta x^{W|E} -x_i) + \nonumber \\
        &&Q_{22} \{+1\} x_i ] \nonumber
\end{eqnarray}
%\end{align}
%\end{subequations}
Summing up the above fluxes gives the discretization scheme for the coarse mesh.
We note that the simple five point stencil used on the fine mesh extends to a 9 point stencil on the coarse mesh.
In practice the proper handling of the sum on the padded domain requires some integer arithmetics.
Note that the above summation leads to a nine point stencil.

\subsection{\textcolor{black}{Multigrid strategy}}
\begin{figure}
\begin{center}
\tikzstyle{block} = [rectangle, draw,fill=gray!20,text width=1.4em,
    text centered, rounded corners, minimum height=1.8em,thick]
\tikzstyle{block2} = [rectangle, draw,color=gray,fill=gray!20,text width=1.4em,
    text centered, rounded corners, minimum height=1.8em,thick]

\begin{tikzpicture}[node distance = 10mm, auto]
%drawing boxes
\node [block] (1_h1) {$\Omega^h$};
\node [block, below right of=1_h1, node distance = 13mm] (1_h2) {$\Omega^{2h}$};
\node [block, below right of=1_h2, node distance = 13mm] (1_h4) {$\Omega^{4h}$};
\node [block, below right of=1_h4, node distance = 13mm] (1_h8) {$\Omega^{8h}$};

\node [block, above right of=1_h8, node distance = 13mm] (2_h4) {$\Omega^{4h}$};
\node [block, above right of=2_h4, node distance = 13mm] (2_h2) {$\Omega^{2h}$};
\node [block, above right of=2_h2, node distance = 13mm] (2_h) {$\Omega^{h}$};

\draw [thick,orange,->] (1_h1) |- (1_h2);
\draw [thick,orange,->] (1_h2) |- (1_h4);
\draw [thick,orange,->] (1_h4) |- (1_h8);

\draw [thick,black,->] (1_h8) -| (2_h4);
\draw [thick,black,->] (2_h4) -| (2_h2);
\draw [thick,black,->] (2_h2) -| (2_h);

\node [left of=1_h1] {(a)};

\end{tikzpicture}

\begin{tikzpicture}[node distance = 10mm, auto]
%drawing boxes
\node [block] (1_h1) {$\Omega^h$};
\node [block2, below right of=1_h1, node distance = 13mm] (1_h2) {$\Omega^{2h}$};
\node [block2, below right of=1_h2, node distance = 13mm] (1_h4) {$\Omega^{4h}$};
\node [block, below right of=1_h4, node distance = 13mm] (1_h8) {$\Omega^{8h}$};

\node [block2, above right of=1_h8, node distance = 13mm] (2_h4) {$\Omega^{4h}$};
\node [block2, above right of=2_h4, node distance = 13mm] (2_h2) {$\Omega^{2h}$};
\node [block, above right of=2_h2, node distance = 13mm] (2_h) {$\Omega^{h}$};

\draw [thick,orange,->] (1_h1) |- (1_h8);
\draw [thick,black,->] (1_h8) -| (2_h);

\node [left of=1_h1] {(b)};

\end{tikzpicture}
\end{center}
\caption{\label{fig:sparse_cycle}
Schematic diagram of multigrid cycles: $\Omega^{xh}$ denote discretizations of different levels, orange, and black arrows correspond to restriction and prolongation steps, respectively.
(a) V cycle using the conventional 2h restriction and prolongation steps.
(b) sparse V cycle, when $\Omega^{2h}$, and $\Omega^{4h}$ discretizations are bypassed.
}
\end{figure}
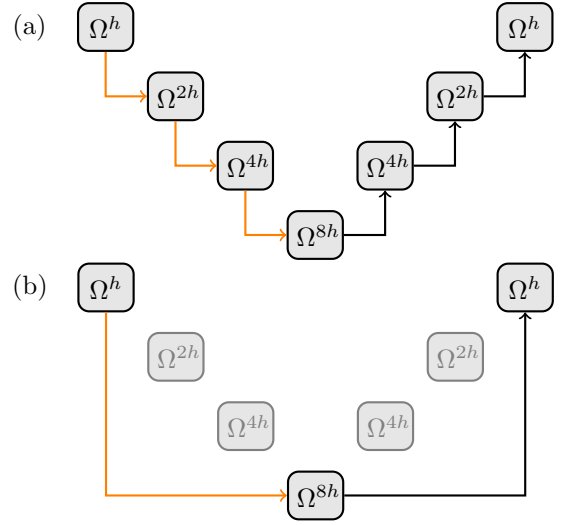

\textcolor{black}{
Multigrid strategy including cycling (i.e. V, W or F-cycling) and the coarsening strategy is a subtle issue, and often examined in term of the arithmetic complexity \cite{Zhao_2005_NMPDE,Brenner_2004_MC}.
The effect of cycling strategy is also known to have effect on parallel performance of the scheme \cite{Lesly_1994_TREP}, however the coarsening strategy is rarely examined in this respect.
It is known, that not all discretization levels are equally important in reducing residual errors, the optimal coarsening may not even uniform\cite{Piquet_2000_NA}.
Herein we wish to investigate the effect of a more aggressive coarsening strategy, in terms of the computational complexity and the parallel scalability in the latency limit.
}
Fig. \ref{fig:sparse_cycle} shows a comparison of a common $2h$ V cycle and its two level "sparse" counterpart.
Note that in the latter case, the number of restriction and prolongation steps are greatly reduced.
In the next section we are going to investigate the effect of substituting common V cycles with sparse ones.

\section{Implementation and benchmarking}

\textcolor{black}{
Our ACM solver was implemented using NVIDIA Cuda and supports only a single device.
Compute tasks for all grid levels run largely on the graphics card, except when global reduction is needed (i.e. convergence check).
We use the five point Laplacian stencil on the fine grid, and the coarser grid discrete equations are generated by summing fine grid equations over the coarser grid, that also gives a 5 point stencil.
Iterations are executed parallel by using a red-black Gauss-Seidel (GS) method.
Reductions are executed in two steps starting on graphics board, and partially reduced data are transferred to the main memory to finish reduction.
} 

\textcolor{black}{In our improved ISMG solver all computations on the fine grid are implemented using the OpenCL streaming language and executed on accelerator cards accordingly.
Our ACM implementation may also use multiple graphics cards over MPI protocol.
While the fine grid equations are solved on the graphics boards, coarse grid equations are implemented in standard C language, and run on the host machine(s).
This adds less latency when iterating on the small sized coarse grid, and avoids poor utilization of the graphics card.
We note, that the free time-slot on the GPU can be utilized for arithmetically intensive pre or post-processing, or for solving coupled scalar equations (i.e. convenction-diffusion equation) that is common in many fluid dynamics applications.
As for the ACM solver we use parallel red-black Gauss-Seidel iterations on the fine grid, but we use a block iterative GS method \cite{Allard_2009_HPCC} on coarse grid level.
GS iteration blocks are aligned with the MPI parallel layout, and simplifies to the serial GS when running on a single node.
The coarse grid solver was implemented as host code, using C language.
We note however, that using parallel multicolor iteration schemes \cite{Adams_1988_SIAM_JNA}, are also an option, which may increase computing efficiency of the CPU device on the expense of extra synchronization, and may increase latency.
}

The OpenCL language allows an efficient use of SIMD engines via defining blocks of operations (work-groups) that are executed asynchronously.
Work-groups can be either one, two or three-dimensional, thus explicit finite difference stencils are straightforward to implement; work-groups are paired with blocks of gridded data (data "tiles" in 2D) \cite{MICIKEVICIUS_2009_PROC}.
\textcolor{black}{However, we have to note that the implementation in 2D and 3D may differ because of the limited number of registers, and out of order execution is not commonly available on graphics hardware.
For the sake of simplicity and easier programming we have chosen a 2D implementation.
Despite our results are in 2D, the extrapolation to estimate 3D performance can be done since convergence rate for GS iterations are depend on the linear size of the simulation domain.
}
In our software code each workgroup corresponds to a coarse-grid cell, thus code complexity is significantly reduced.
The best computing efficiency can be achieved, when the tiles/work-groups are fitted to the SIMD width of the compute device.
Using an optimal tile size on GPUs, the amount of data to be processed on the coarse grid is greatly reduced (i.e., by a factor of $256$ for a $16\times16$ tile size).

Code performance and our multigrid scheme was analyzed on various hardware platforms.
We wish to emphasize that we have not put extra efforts in optimizing for specific platforms.
Our aim was to measure the efficiency of the MG scheme on a more generic way, and give estimates on possible bottlenecks.
%We have not run the extra mile to prepare the whole scheme in native C, however some basic tests we have performed suggest that OpenCl  gives at least a factor of two in efficiency, when compared to a "naive" C code using openMP parallelism.

\subsection{OpenCl kerel performance}
Our OpenCl compute kernel implementations use global memory, except reduction kernels (i.e. restrictions steps in MG) \textcolor{black}{which also use local address space}.
This means that the kernels are quasi-optimal for CPUs, and further speedup is possible on graphics processors.
The reason for using global memory access is its simplicity.
We wanted to use exactly the same compute kernels on each platform, and we have not found a generic way to utilize the full potential of the specific hardwares.
Thus, our kernel timings can easily be achieved by a scientist or fluid dynamics expert, who is familiar with GPU programming at a basic level.
%Also the migration from 2D to 3D in the future will be even complicated, since in some kernels cache requirement exceeds hardware limits.
%Although this can be bypassed by using cache efficient algorithms \cite{MICIKEVICIUS_2009_PROC}
Implementations of the finite difference stencils may also differ in flexibility of handling boundary conditions and grid size.
Flexibility of the code often leads to conditional statements, which may dramatically affect code performance.
In practice, conditional statements can mostly be replaced with integer, or bit-wise integer arithmetics, but leads to extra operations.
We also note that architectures and platforms may significantly differ in performing integer operations that may result in differences in computing time for codes with flexible boundary conditions.
Our kernel implementation uses automated padding if the grid size is not a multiple of the tile size, and can hold symmetric, Dirichlet and periodic boundary conditions.
We use wall times of kernel execution as the measure of the performance.
The results benchmarked on some recent hardware are summarized in Table \ref{table:kernel_time}.

\begin{table}[ht]
\caption{single chip kernel execution wall clock time/million grid points vs. work-group size}
\label{table:kernel_time}
\begin{tabular}{cccccccc}
\hline
\hline
i5-2500K		& $4 \times 4$	& $8 \times 8$	& $16\times16$	& $32\times32$ \\
\hline
Navier-Stokes		& 9.2 ms		& 10.1 ms		& 7.4 ms	& 7.9 ms \\
red-black GS		& 4.2 ms		& 3.3 ms		& 3.0 ms	& 3.2 ms \\
divergence		& 2.6 ms		& 1.7 ms		& 1.2 ms	& 1.2 ms \\
restriction 		& 6.1 ms		& 7.8 ms		& 8.7 ms	& 12.4 ms \\
prolongation 		& 8.8 ms		& 7.9 ms		& 7.7 ms	& 7.85 ms \\
%\hline
%\hline
%GTX 590			& $4 \times 4$	& $8 \times 8$	& $16\times16$	& $32\times32$ \\
%\hline
%predictor		& 1.41 ms		& 1.01 ms		& 0.73 ms	& - \\
%pressure		& 0.84 ms		& 0.34 ms		& 0.23 ms	& - \\
%divergence		& 0.66 ms		& 0.30 ms		& 0.20 ms	& - \\
%reduction 		& 0.87 ms		& 0.27 ms		& 0.19 ms	& - \\
%interpolation 		& - ms		& - ms		& - ms	& - ms \\
\hline
\hline
GTX 680			& $4 \times 4$	& $8 \times 8$	& $16\times16$	& $32\times32$ \\
\hline
Navier-Stokes		& 1.18 ms		& 0.68 ms		& 0.38 ms	& 0.39 ms\\
red-black GS		& 0.68 ms		& 0.22 ms		& 0.13 ms	& 0.15 ms\\
divergence		& 0.51 ms		& 0.20 ms		& 0.12 ms	& 0.13 ms\\
restriction 		& 0.99 ms		& 0.30 ms		& 0.17 ms	& 0.27 ms\\
prolongation 		& 1.06 ms		& 0.35 ms		& 0.33 ms	& 0.41 ms\\
%\hline
%\hline
%Radeon 5970		& $4 \times 4$	& $8 \times 8$	& $16\times16$	& $32\times32$ \\
%\hline
%predictor		& 5.83 ms		& 1.55 ms		& 1.52 ms	& - \\
%pressure		& 0.62 ms		& 0.19 ms		& 0.15 ms	& - \\
%divergence		& 0.55 ms		& 0.17 ms		& 0.14 ms	& - \\
%reduction 		& 2.51 ms		& 0.81 ms		& 0.60 ms	& - \\
%interpolation 		& - ms		& - ms		& - ms	& - ms \\
\hline
\hline
Radeon 7970		& $4 \times 4$	& $8 \times 8$	& $16\times16$	& $32\times32$ \\
\hline
Navier-Stokes		& -			& 0.33 ms		& 0.28 ms	& - \\
red-black GS		& -			& 0.07 ms		& 0.06 ms	& - \\
divergence		& -			& 0.09 ms		& 0.06 ms	& - \\
restriction 		& -			& 0.21 ms		& 0.12 ms	& - \\
prolongation 		& -			& 0.31 ms		& 0.26 ms	& - \\
\hline
\end{tabular}
\end{table}

According to the benchmarks, execution times are decreasing as work-group size increases on all hardware and saturates at the SIMD width.
Its variation can be significant on GPUs.
We note that with increasing restriction length we have double benefit; not just the computational cost of solving an iteration on the coarse mesh is decreasing, but also iterating on the fine grid becomes more efficient.

\subsection{\textcolor{black}{Accuracy test}}
\textcolor{black}{
In order to test accuracy and the righteousness of our implementation we have chosen the well-known lid-driven cavity test first published by Ghia \cite{Ghia_1982_JCP}.
Our computations were executed until steady solution was reached.
The extrema of the velocities along the centerlines are in reasonable agreement with the values reported in previous works (see table \ref{table:lid_driven}.).
We also note that deviations among the reported values are within the 32 bit floating point precision (after approx. 40000 non-dimensional time), indicating that our implementation using single precision is reasonable.
}
\begin{table}[ht]
\caption{Extrema of the velocities through the centerlines of the cavity, at Re = 1000}
\label{table:lid_driven}
\begin{tabular}{lccccccc}
\hline
\hline
reference		& grid	& $u_{max}$	& $v_{min}$	& $v_{max}$ \\
\hline
\cite{Ghia_1982_JCP}		& $129\times129$		& 0.3829		& -0.5155	& 0.3710 \\
\cite{Deng_1994_CAF}	stagg.	& $128\times128$		& 0.3805		& -0.5173	& 0.3688 \\
\cite{Bruneau_1990_JCP}		& $256\times256$		& 0.3764		& -0.5208	& 0.3665 \\
\cite{Vanka_1986_JCP}		& $321\times321$		& 0.3870		& -	& - \\
\cite{Botella_1998_CAF} spectral	& $160\times160$		& 0.3886		& -0.5271	& 0.3769 \\
present work	& $512\times512$		& 0.3781		& -0.5142	& 0.3659 \\
\hline
\end{tabular}
\end{table}

\subsection{Multigrid performance}

We have measured the computational cost (arithmetic intensity) of the MG scheme using the average number of computed Laplacians corresponding to a fine grid control volume in a time increment ($N_{Lap}$).
This can be used to compare different MG schemes in an ideal case, where the network communication and the memory bandwidth are not limiting the calculation.

Although the computational cost of performing an iteration is an essential feature in evaluating a numerical method, the number of iterations required to reach convergence criteria is more fundamental, when evaluating scalability on large distributed memory computers.
While the raw computational power can be increased by stacking computers, the interconnect latency and bandwidth have strict limitations.
When performing finite difference iterations, the communication can be hidden up to the point, when compute time decreases to the network communication time.
Down to this limit, the compute power can be fully utilized: theoretically the speedup is proportional to the number of processors.
When pushing to this limit (i.e. by increasing the number of the compute devices), the total compute time can be reduced only via decreasing the iteration count, if we maintain computing efficiency.
We  use the total GS iteration count ($I_{t}$) and the related total network communication count (synchronizations) per time step ($NCC_t$) as the measure of performance on distributed memory clusters.
We distinguish iteration numbers on fine ($I_f$) and coarse grids ($I_c$) since they behave differently with varying numerical parameters,
and network communication characteristics also varies with grid size.
The actual choice for the solver also influences network communication requirements: the simple \textcolor{black}{block iterative} GS involves 1 synchronization per iteration that doubles when red-black reordering is applied.
%Coarse and fine grid iterations have different communication characteristics, since on coarse mesh less data is needed to be sent.
Besides, when compute devices are external cards, an extra, compute device to host transfer is also involved.

\subsubsection{Benchmark setup I.: Shear driven laminar flow in cavity}
\begin{figure}[ht!]
\includegraphics[width=1\linewidth]{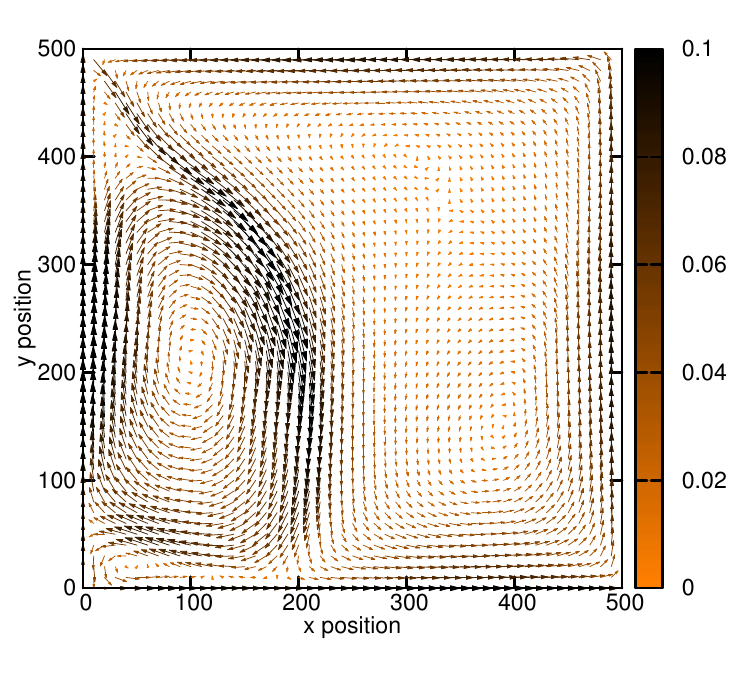}
\caption{\textcolor{black}{Illustrative vector plot for} benchmark I.: counter clockwise shear driven cavity flow with reversed velocity on the left, kinematic viscosity $\eta=0.1$, velocity at the boundaries $v_0=0.1$, grid spacing $h=1$, time-step $\Delta t=1$ grid size $500\times500$.
Color code shows the magnitude of the velocity vectors.
}\label{fig:shear_driven_cavity_flow}
\end{figure}
In this section we use laminar flow in a rectangular cavity as a benchmark.
\textcolor{black}{To make tests more challenging, we put a twist into the commonly used lid-driven cavity by prescribing opposite velocities on the fixed walls.}
The setup and the numerical parameters are shown in Fig.\ref{fig:shear_driven_cavity_flow}.
When using advanced multigrid methods, often complicated patterns are used to restrict and prolongate among various grid-levels, such as v-cycles and w-cycles.
Here, the simplest accommodative cycles strategy is used: %, which decide adaptively on changing the grid:
Iteration is started on the coarsest grid, and when the absolute maximum of the residual is decreased below a limit, the solution is prolongated to the finest grid.
We can vary this criterion to find a quick solution.
Then, iteration continues on the fine grid until reaching the fine grid convergence criterion.
In our benchmark we require the residual to decrease below $10^{-6}$ on the fine grid, and vary the coarse grid criterion.
Reaching this criterion requires significant efforts, but we are still comfortably far from the limits of the 32 bit float arithmetics.
\begin{figure}[ht!]
	\begin{center}
	\includegraphics[width=1\linewidth]{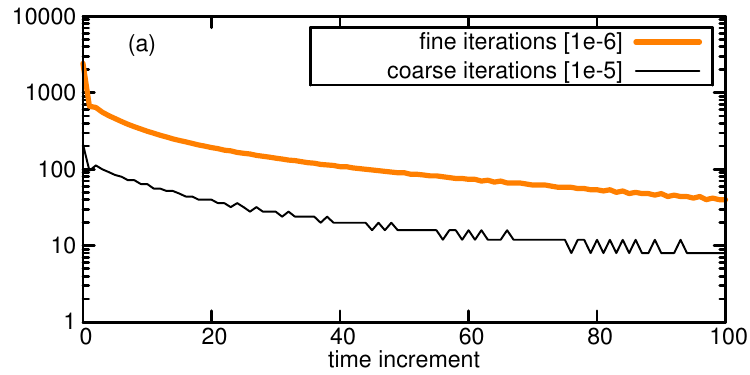}
	\end{center}
	\begin{center}
	\includegraphics[width=1\linewidth]{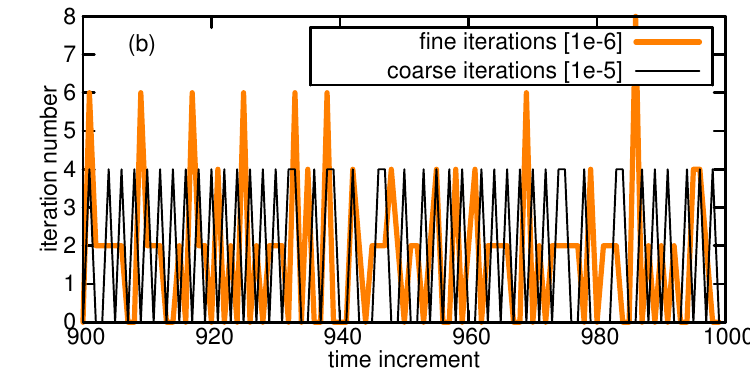}
	\end{center}
\caption{Benchmark I.: iteration numbers on the coarse and the fine grid vs. time increment. Convergence criteria was set to $10^{-6}$ and $10^{-5}$ for the fine and the coarse grid iterations, respectively.
(a) very strong transient in first 100 time increments, (b) moderate transient between 900-1000 increments.
}\label{fig:iter_vs_time}
\end{figure}

We have compared our results to the ACM scheme \cite{Guerrero_2000_thesis} that is using the conventional 2h restriction length.
On the coarsest mesh, we have applied the same convergence criterion as for the ISMG scheme.
We have performed heuristic tests to find a quasi-optimal cycle: V cycles without pre-smoothing, and one post-smoothing per grid level was found to be a good candidate.
The grid depth was varied, thus we have the ACM substitutes for our ISMG scheme using varying restriction length (e.g. tile size).
In the case of our ACM implementation we denote only the coarsest mesh iterations as "coarse", and all the others are accounted as "fine".
Iteration counts at each grid level can be calculated, since the number of the smoothing steps are the same, except for the coarse mesh.

We have measured scalability ($NCC$s) and \textcolor{black}{arithmetic intensity} ($N_{Lap}$) that are summarized in table \ref{table:laminar_500_ISMG} for the  ISMG scheme and in Table \ref{table:laminar_500_ACM} for the ACM scheme, respectively.
The optimum values we have found are highlighted with bold characters.
We have measured similar computational cost while our ISMG method requires six times less synchronizations.
We note that the greatly reduced synchronization count also results in less kernel launches that is preferable in streaming computer languages.
Besides, restriction and prolongation kernels are executed fewer times, however it is difficult to make a direct comparison on this basis, since it depends largely on the multigrid cycling pattern.

The low iteration counts we measured are partly due to the fact that in the ISMG scheme the coarse grid solution approximates the fine grid discrete solution.
Fig. \ref{fig:iter_vs_time}.b shows that in the moderate transient regime after prolongating to the fine grid the solution occasionally passing convergence criteria with no iterations.
The same can be observed the way back, sometimes no coarse iteration is required after a converged solution of the previous time-step is applied in the predictor step.

\subsubsection{Benchmark setup II.: Turbulent jet}
\begin{figure}[ht!]
\includegraphics[width=1\linewidth]{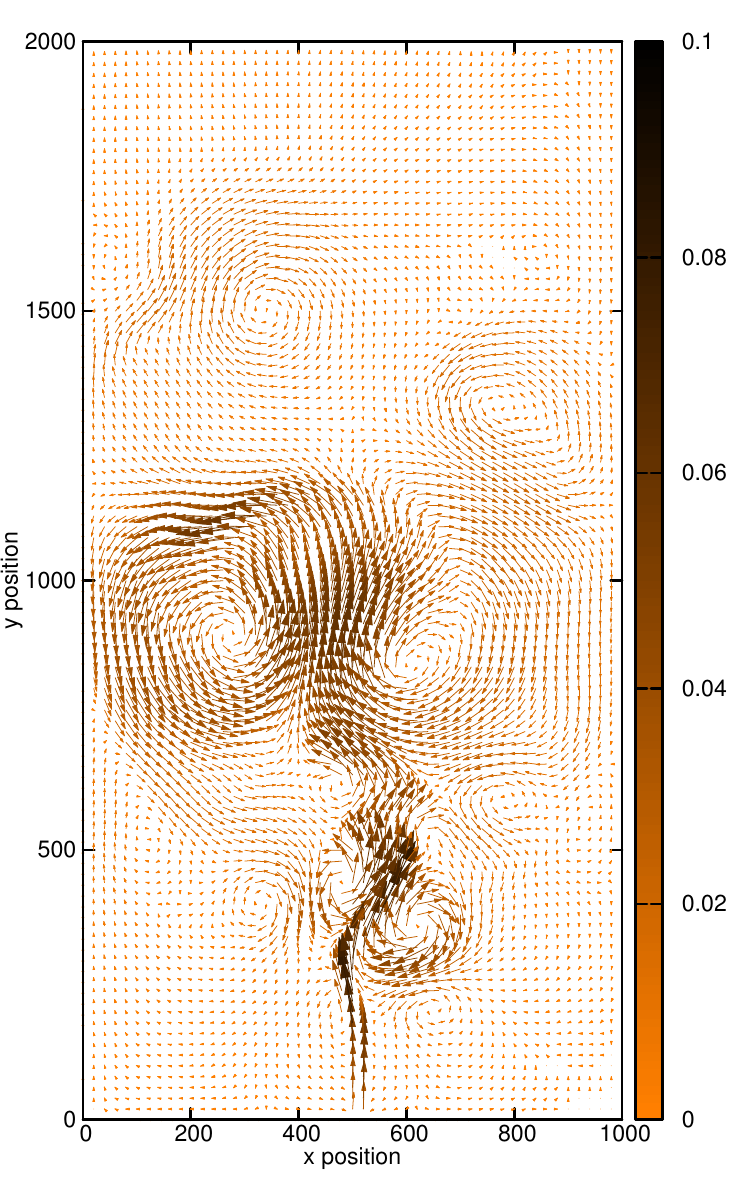}
\caption{Benchmark II.: 2D turbulent jet flow in a rectangular cavity,
no slip walls at the left and right boundaries;
symmetry for velocity, and fixed pressure at the top boundary,
and a tiny, 16 control volume wide inlet with $v_0=0.1$ is defined in the middle of the no slip bottom boundary,
kinematic viscosity $\eta=0.01$,
grid spacing $h=1$, time-step $\Delta t=1$ grid size $1000\times2000$.
Color code shows the magnitude of the velocity vectors.
}\label{fig:tubulent_jet_flow}
\end{figure}
The previous test setup can be regarded as a synthetic benchmark, since the pressure we calculate is deterministic, therefore, one can tune the setup to provide unbiased benchmarks for the actual scheme (e.g. by varying the grid-size).
In this section, we present benchmarks for Direct Numerical Simulation (DNS) of a turbulent jet.
The actual setup, and the numerical parameters applied are shown in Fig. \ref{fig:tubulent_jet_flow}.
Due to the chaotic nature of turbulence, it would be very difficult to tweak the test setup to show unbiased benchmarks, it can be accounted as a real-life application of our scheme.
We have increased the grid size to $1000\times2000$.
Here we have maximized the number of GS iterations at 20000 for both ACM and ISMG computations to prevent convergence stagnation, however we believe this can be bypassed by using a more advanced cycling strategy.
The linear sizes we have used in this setup are close to the supercomputing scale 3D direct numerical simulations.

Due to the increase in the linear size the small grid depth computations for ACM scheme and the corresponding $4\times4$ and $8\times8$ tile sizes in the ISMG scheme the solution is far from optimal and requires excessive compute time, therefore results are not presented for these cases.
The results for 5-6 grid levels of the ACM and the corresponding $16\times16$ and $32\times32$ tile size ISMG schemes are presented in Table \ref{table:tubulent_jet}.
\textcolor{black}{For comparison we also have made a benchmark for this setup using an aggressive coarsening strategy, but with the same five point coarse grid stencil as is usual in geometric multigrid (GMG) schemes.
Results (last two columns of Table \ref{table:tubulent_jet}) provide insight how arithmetic intensity decreases when using our interpolated stencil methodology.
}
In the case of the ACM scheme the optimum has been found at 6 grid levels contrary to the ISMG scheme, where the optimum has been found at $16\times16$ tile size (corresponds to 5 grid levels of the ACM method).
Comparing the optima of the two methods, the ISMG scheme maintained the sixfold advantage in scalability ($NCC_t$), but at the cost of approximately twofold increase in the computational need ($N_{Lap}$).
We note that in case of the applied cycling strategy this computational overhead is compensated by reduced restriction and prolongation counts, and their more efficient execution.
\textcolor{black}{When comparing the interpolated stencil to the conventional five point Laplacian, we have measured $26\%$ reduction in the arithmetic intensity of the iterations.}

\subsubsection{Benchmark setup III.: Strong scaling}
Finally, we asses multi-GPU performance and efficiency of our algorithm for different problem sizes.
The benchmark setup is a channel flow having equally spaced jets from a side.
Each jets having width of $25 h$, (just as for benchmark II), and their spacing is $200 h$.
Above construction allows a rough comparison of different problem sizes.
However the solution is inherently stochastic, averaging the first $20000$ time-step gives a reasonable basis for the comparison.

Our tests has been performed on two nodes totaling eight graphics processors, and connected with Mellanox SDR infinband.
The detailed configuration for a nodes is:

CPU: Dual Intel Xeon E5530 (2.4 GHz), 2 x Nvidia GTX 590 boards (4 x GF110 GPU per node).
Our benchmarks were run under CentOS 6.5, and using Cuda 5.5 along with Nvidia driver version 319.37.

\begin{figure}[ht!]
\includegraphics[width=1\linewidth]{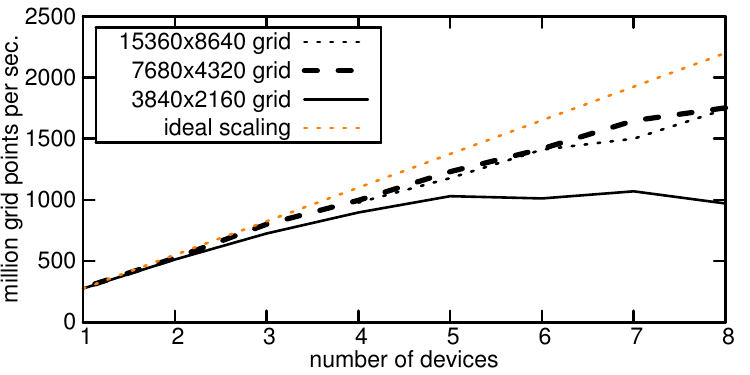}
\caption{strong scaling:
The average computing speed for a time-step for different problem sizes.
The setup consists of equally spaced jets, similar to that of soap film experiments.
}\label{fig:strong_scaling}
\end{figure}

The benchmarks are summarized in Fig. \ref{fig:strong_scaling}.
The largest benchmark (15360x8640) saturates the 1.5 GB memory per GPU, therefore strong scaling curves showe for 4-8 graphics cores.
For large enough setups, we have measured up to $80\%$ efficiency at 8 graphics cores.
Losing efficiency comes partly from gathering communication when checking convergence, and partly from transferring the coarse grid problem forth and back between the GPU memory and the main memory.
The latter can be hidden, if some coupled phenomenon (i.e. thermal convection and diffusion) are needed to solve along with the fluid flow.
We also note, that in three dimensions, the efficiency curve is expected to be even shallower due to the more rapid decrease in the problem size to coarse meshes.

\section{Concluding remarks}
We have presented an interpolated stencil multigrid method that is efficient when 16-32h wide restriction step is applied.
The stencil allows the replacement of a deep grid hierarchy (optimally 5-6 grid levels) with a two grid scheme, while keeping the number of iterations low.
The scheme was benchmarked against an ACM scheme with deep grid hierarchy.
It is shown that ISMG needs significantly less iterations compared to a conventional ACM scheme, since it involves less synchronizations (i.e. kernel calls and network communications).
Also, at larger grid sizes the raw computational cost (the number of the computed Laplacians) is increased by a factor of two, however it is overcompensated by the reduced number of prolongation and restriction kernel calls.
We have measured the parallel efficiency near $80\%$ for a cluster includes 8 graphics cores.
For further fine tuning our scheme, we plan to replace the simple one-step accommodative cycling strategy with more complex cycling patterns to further decrease computational cost.
Using our method with higher order interpolation is also an issue that we wish to investigate.
Since Gauss-Seidel has convergence rate of $O(N^2)$, where $N$ is the linear size of the grid, we expect even faster convergence in three dimensional cases.
Besides, in 3D the coarse graining results in faster reduction of the coarse problem size, thus helping the more effective use of the parallel hardware and better strong scaling.

\section{Acknowledgements}
We thank L. Gr\'an\'asy for illuminating discussions and for supporting our research.
This work has been supported by the EU FP7 Collaborative Project ‘‘EXOMET’’ (contract no. NMP-LA-2012-280421, co-funded by ESA), and by the ESA MAP project ‘‘GRADECET’’ (Contract No. 4000104330/11/NL/KML).
G. Tegze is a grantee of the J\'anos Bolyai Scholarship of the Hungarian Academy of Sciences.
We also thank the MTA Wigner RCP's GPU Lab for for supporting us with hardware and programming advices.

\appendix

\section{Benchmarks}

\begin{table}[ht]
\caption{Benchmark I: the network communication counts in the fine ($LCSS_f$) and the coarse grid ($LCSS_c$), their sum ($LCSS_t$) and the raw computational cost ($N_{Lap}$) in the first 1000 time increment using the ISMG scheme for various tile sizes and coarse grid convergence criterion.
Quasi optimal values are highlighted using bold characters.}
\label{table:laminar_500_ISMG}
\begin{tabular}{crrrrrrr}
coarsest		& $4h$	& $8h$	& $16h$	& $32h$ \\
\hline
\hline
$NCC_{f}$	\\
\hline
$1\times 10^{-6}$	&	5.0	&	9.7	&	21.4	&	50.8	\\
$1\times 10^{-5}$	&	45.2	&	10.7	&	22.4	&	51.2	\\
$2\times 10^{-5}$	&	166.1	&	13.4	&	23.8	&	51.6	\\
$3\times 10^{-5}$	&	238.1	&	57.0	&	26.2	&	52.0	\\
$4\times 10^{-5}$	&	262.2	&	78.3	&	28.9	&	52.6	\\
$5\times 10^{-5}$	&	273.4	&	97.3	&	32.6	&	53.0	\\
$1\times 10^{-4}$	&	296.9	&	221.0	&	85.4	&	59.5	\\
\hline
\hline
$NCC_{c}$	\\
\hline
$1\times 10^{-6}$	&	72.8	&	47.4	&	31.2	&	17.4	\\
$1\times 10^{-5}$	&	11.2	&	10.1	&	7.7	&	5.1	\\
$2\times 10^{-5}$	&	4.3	&	6.5	&	5.2	&	3.8	\\
$3\times 10^{-5}$	&	2.2	&	8.5	&	4.2	&	3.2	\\
$4\times 10^{-5}$	&	1.5	&	2.5	&	3.6	&	2.9	\\
$5\times 10^{-5}$	&	1.1	&	1.9	&	3.3	&	2.8	\\
$1\times 10^{-4}$	&	0.5	&	0.6	&	3.1	&	2.3	\\
\hline
\hline
$NCC_t$ \\
\hline
$1\times 10^{-6}$	&	77.8	&	57.1	&	52.6	&	68.2	\\
$1\times 10^{-5}$	&	56.4	&	{\bf20.8}	&	30.0	&	56.3	\\
$2\times 10^{-5}$	&	170.3	&	{\bf19.8}	&	29.0	&	55.4	\\
$3\times 10^{-5}$	&	240.2	&	65.5	&	30.4	&	55.3	\\
$4\times 10^{-5}$	&	263.7	&	80.7	&	32.5	&	55.5	\\
$5\times 10^{-5}$	&	274.6	&	99.2	&	35.9	&	55.8	\\
$1\times 10^{-4}$	&	297.4	&	221.6	&	88.5	&	61.8	\\
\hline
\hline
$N_{Lap}$ \\
\hline
$1\times 10^{-6}$	&	7.0	&	{\bf5.6}	&	10.8	&	25.4	\\
$1\times 10^{-5}$	&	23.3	&	{\bf5.5}	&	11.2	&	25.6	\\
$2\times 10^{-5}$	&	83.3	&	6.8	&	11.9	&	25.8	\\
$3\times 10^{-5}$	&	119.2	&	28.6	&	13.1	&	26.0	\\
$4\times 10^{-5}$	&	131.2	&	39.2	&	14.5	&	26.3	\\
$5\times 10^{-5}$	&	136.8	&	48.7	&	16.3	&	26.5	\\
$1\times 10^{-4}$	&	148.5	&	110.5	&	42.7	&	29.8	\\
\hline
\end{tabular}
\end{table}

\begin{table}[ht]
\caption{Benchmark I: the network communication counts on the fine ($LCSS_f$) and the coarse grid ($LCSS_c$), their sum ($LCSS_t$) and the raw computational cost ($N_{Lap}$) in the first 1000 time increment using the ACM scheme for various grid-depths and coarse grid convergence criterion.
Quasi optimal values are highlighted using bold characters.}
\label{table:laminar_500_ACM}
\begin{tabular}{crrrrrrr}
coarsest		& 4h	& 8h	& 16h	& 32h \\
\hline
\hline
$NCC_{f}$\\
\hline
$1\times 10^{-6}$	&	6.6	&	14.2	&	24.8	&	170.8	\\
$1\times 10^{-5}$	&	48.2	&	15.0	&	25.4	&	170.8	\\
$2\times 10^{-5}$	&	182.2	&	18.2	&	26.0	&	170.8	\\
$3\times 10^{-5}$	&	259.8	&	51.2	&	27.2	&	170.8	\\
$4\times 10^{-5}$	&	294.2	&	81.4	&	28.8	&	170.8	\\
$5\times 10^{-5}$	&	314.0	&	115.4	&	30.8	&	170.8	\\
$1\times 10^{-4}$	&	353.0	&	216.6	&	70.8	&	170.8	\\
\hline
\hline
$NCC_{c}$\\
\hline
$1\times 10^{-6}$	&	1012.0	&	1177.0	&	1151.4	&	0.0	\\
$1\times 10^{-5}$	&	177.8	&	273.0	&	326.6	&	0.0	\\
$2\times 10^{-5}$	&	77.4	&	165.8	&	211.2	&	0.0	\\
$3\times 10^{-5}$	&	46.2	&	110.8	&	161.4	&	0.0	\\
$4\times 10^{-5}$	&	33.0	&	81.4	&	132.4	&	0.0	\\
$5\times 10^{-5}$	&	25.6	&	62.4	&	112.6	&	0.0	\\
$1\times 10^{-4}$	&	11.8	&	26.2	&	57.2	&	0.0	\\
\hline
\hline
$NCC_t$ \\
\hline
$1\times 10^{-6}$	&	1018.6	&	1191.2	&	1176.2	&	170.8	\\
$1\times 10^{-5}$	&	226.0	&	288.0	&	352.0	&	170.8	\\
$2\times 10^{-5}$	&	259.6	&	184.0	&	237.4	&	170.8	\\
$3\times 10^{-5}$	&	306.0	&	161.8	&	188.6	&	170.8	\\
$4\times 10^{-5}$	&	327.2	&	162.8	&	161.0	&	170.8	\\
$5\times 10^{-5}$	&	339.6	&	177.8	&	143.6	&	170.8	\\
$1\times 10^{-4}$	&	364.8	&	242.8	&	{\bf 128.0}	&	170.8	\\
\hline
\hline
$N_{Lap}$ \\
\hline
$1\times 10^{-6}$	&	33.3	&	12.2	&	6.3	&	22.7	\\
$1\times 10^{-5}$	&	17.6	&	{\bf 5.3}	&	{\bf 4.8}	&	22.7	\\
$2\times 10^{-5}$	&	48.0	&	{\bf 5.1}	&	{\bf 4.7}	&	22.7	\\
$3\times 10^{-5}$	&	66.4	&	11.5	&	{\bf 4.8}	&	22.7	\\
$4\times 10^{-5}$	&	74.6	&	17.6	&	{\bf 5.0}	&	22.7	\\
$5\times 10^{-5}$	&	79.3	&	24.5	&	{\bf 5.3}	&	22.7	\\
$1\times 10^{-4}$	&	88.6	&	45.3	&	11.7	&	22.7	\\
\hline
\end{tabular}
\end{table}

\begin{table}[ht]
\caption{Benchmark II: the network communication counts the fine ($LCSS_f$) and the coarse grid ($LCSS_c$), their sum ($LCSS_t$) and the raw computational cost ($N_{Lap}$) in the first 100000 time increment using the ACM and ISMG schemes for various grid-depths/tile sizes and coarse grid convergence critera.
Quasi optimal values are highlighted using bold characters.}
\label{table:tubulent_jet}
\begin{tabular}{crrrrrrr}
coarsest		& $16h$	& $32h$	& $16h$	& $32h$ & $16h$ & $32h$\\
\hline
\hline
$NCC_{f}$		&	ISMG &	ISMG	&	ACM	&	ACM	&	GMG	&	GMG\\
\hline
$1\times 10^{-6}$	&	3.4	&	9.5	&	14.2	&	46.2	&	3.8	&	9.7\\
$1\times 10^{-5}$	&	3.5	&	8.6	&	4.4	&	138.8	&	4.5	&	7.8\\
$2\times 10^{-5}$	&	-	&	9.5	&	4.4	&	5.4	&	-	&	8.3\\
$3\times 10^{-5}$	&	-	&	8.1	&	4.4	&	5.4	&	-	&	8.0\\
$4\times 10^{-5}$	&	-	&	8.0	&	4.4	&	5.4	&	-	&	10.4\\
$5\times 10^{-5}$	&	-	&	8.7	&	4.6	&	5.4	&	-	&	12.6\\
$1\times 10^{-4}$	&	-	&	13.5	&	-	&	5.4	&	-	&	14.9\\
\hline
\hline
$NCC_{c}$	\\
\hline
$1\times 10^{-6}$	&	5.0	&	7.9	&	428.0	&	494.2	&	1128.4	&	15.1\\
$1\times 10^{-5}$	&	755.4	&	4.1	&	184.0	&	269.2	&	4.0	&	4.8	\\
$2\times 10^{-5}$	&	-	&	3.6	&	62.0	&	67.2	&	-	&	8.1	\\
$3\times 10^{-5}$	&	-	&	2.7	&	198.8	&	121.0	&	-	&	3.1	\\
$4\times 10^{-5}$	&	-	&	2.4	&	240.4	&	125.8	&	-	&	2.6	\\
$5\times 10^{-5}$	&	-	&	2.5	&	272.8	&	104.4	&	-	&	2.8	\\
$1\times 10^{-4}$	&	-	&	2.1	&	-	&	40.4	&	-	&	2.4	\\
\hline
\hline
$NCC_t$ \\
\hline
$1\times 10^{-6}$	&{\bf8.4}	&	17.4	&	442.4	&	540.4	&	1132.1	&	24.8\\
$1\times 10^{-5}$	&	758.9	&	12.8	&	188.2	&	408.0	&	{\bf8.5}	&	12.6\\
$2\times 10^{-5}$	&	-	&	13.1	&	66.4	&	72.6	&	-	&	16.4\\
$3\times 10^{-5}$	&	-	&{\bf10.7}	&	203.0	&	126.4	&	-	&	11.0\\
$4\times 10^{-5}$	&	-	&{\bf10.4}	&	245.0	&	131.2	&	-	&	13.0\\
$5\times 10^{-5}$	&	-	&	11.2	&	277.2	&	109.8	&	-	&	15.4\\
$1\times 10^{-4}$	&	-	&	15.6	&	-	&	{\bf45.8}	&	-	&	17.3\\
\hline
\hline
$N_{Lap}$ \\
\hline
$1\times 10^{-6}$	&{\bf1.7}	&	4.8	&	3.2	&	6.4	&	6.3	&	4.9\\
$1\times 10^{-5}$	&	4.7	&	4.3	&	1.1	&	18.6		&	{\bf2.3}	&	3.9\\
$2\times 10^{-5}$	&	-	&	4.8	&	{\bf0.8}	&	{\bf0.8}	&	-	&	4.2\\
$3\times 10^{-5}$	&	-	&	4.1	&	1.1	&	{\bf0.8}	&	-	&	4.0\\
$4\times 10^{-5}$	&	-	&	4.0	&	1.2	&	{\bf0.8}	&	-	&	5.2\\
$5\times 10^{-5}$	&	-	&	4.4	&	1.3	&	{\bf0.8}	&	-	&	6.3\\
$1\times 10^{-4}$	&	-	&	6.8	&	-	&	{\bf0.7}	&	-	&	7.5\\
\hline
\end{tabular}
\end{table}

%% References


\bibliographystyle{elsarticle-num}
\bibliography{references}






% \end{thebibliography}

\end{document}